# Strategies for Feature-Assisted Development of Topology Agnostic Planar Antennas Using Variable-Fidelity Models

Adrian Bekasiewicz[1], Khadijeh Askaripour[1], Mariusz Dzwonkowski[1,2], Tom Dhaene[3], and Ivo Couckuyt[3]

[1] Faculty of Electronics, Telecommunications and Informatics, Gdansk University of Technology, Gdansk, Poland, adrian.bekasiewicz@pg.edu.pl
[2] Department of Radiology Informatics and Statistics, Faculty of Health Sciences, Medical University of Gdansk, Tuwima 15, 80-210 Gdansk, Poland.
[3] Department of Information Technology (INTEC), IDLab, Ghent University-imec, iGent, Technologiepark-Zwijnaarde 126, 9052 Ghent, Belgium

**Abstract**

Design of antennas for contemporary applications presents a complex challenge that integrates cognitive-driven topology development with the meticulous adjustment of parameters through rigorous numerical optimization. Nevertheless, the process can be streamlined by emphasizing the automatic determination of structure geometry, potentially reducing the reliance on traditional methods that heavily rely on engineering insight in the course of antenna development. In this work, which is an extension of our conference paper [1], a specification-oriented design of topologically agnostic antennas is considered by utilizing two strategies focused on bandwidth-specific design and bandwidth-enhanced optimization. The process is embedded within a variable-fidelity framework, where the low-fidelity optimization involves classification of randomly generated topologies, followed by their local tuning using a trust-region algorithm applied to a feature-based representation of structure response. The final result is then tuned using just a handful of high-fidelity simulations. The strategies under consideration were verified on a case study basis concerning automatic generation of three radiators with varying parameters. Benchmarks of the algorithm against more standard optimization methods, as well as comparisons of the obtained topologies with respect to state-of-the-art antennas from literature have also been considered.

## 1. Introduction

Design of modern antennas is an inherently cognitive process. It involves experience-driven development of topology followed by its tuning so as to fulfill the desired performance requirements [2], [3], [4], [41], [42]. When combined with robust optimization, this engineer-in-a-loop approach proved to be useful for the design of new (often unconventional) topologies [5], [6]. Although restricting the antenna development to a specific shape is considered pivotal for ensuring feasibility of the design process, it is also a subject to engineering bias which tends to limit the potential in terms of achieving solutions characterized by unique (perhaps not expected) performance characteristics [5]. These might include, for instance, broadband, or multiband behavior, but also improved radiation capabilities and/or small dimensions [4]-[6]. Consequently, a streamlined design that relies only on numerical optimization methods represents an interesting alternative to the standard antenna development techniques [7], [43], [44].

Automatic antenna generation governed by numerical optimization is a challenging task. From the perspective of geometry, the design can be represented as a set of points (interconnected using, e.g., line-sections, or splines), or in the form of a binary matrix that defines the configuration of primitives (e.g., rectangles) constituting the antenna [7]-[12]. Point-based approaches are capable of supporting geometrically complex topologies. Owing to continuous and free-form nature of coordinates, "evolution" of geometry can be governed by standard numerical optimization algorithms [8], [10]. However, a large number of points and constraints on their distribution are required to ensure flexibility of topology and its consistency (here understood as lack of self-intersections between the coordinate-based curves) [8], [9]. Matrix-based methods naturally comply with limitations by representing topologies as compositions of partially overlapping primitives [7], [11]. Their main bottleneck is a large number of variables required even for relatively simple topologies. Besides activation/deactivation of primitives (based on contents of the matrix), their dimensions also need to be scaled which reveals a mixed-integer nature of the problem [11]. Finally, high dimensionality and the need to evaluate performance based on the expensive electromagnetic (EM) simulations make the discussed universal representations impractical for optimization when conventional algorithms—that require hundreds or even thousands of design evaluations to converge—are considered [5], [10], [11]. Therefore, formulating strategies for identification of geometries that are suitable for optimization remains a significant challenge in terms of both considered approaches.

Population-based algorithms, such as metaheuristics, can potentially help to mitigate the issue pertinent to topology generation by exploring diverse geometries, as they evaluate multiple candidate designs per iteration [13], [14]. However, the computational cost of assessing their responses through EM simulations is tremendous, restricting the practical application of conventional metaheuristics to relatively simple structures represented using a limited number of independent parameters [13], [15]. The problem pertinent to high cost of optimization can be mitigated (to some extent) using surrogate-assisted methods. The latter ones involve substitution of expensive EM models with their cheaper representations [16], [17]. While metaheuristics are useful for generating approximate solutions, their effectiveness diminishes with increasing dimensionality of the problem, often restricting their usefulness to scenarios characterized by a limited number of input parameters [18]. Although extending surrogate model applicability to over a dozen variables is possible, it requires confining the search space to a predefined region of interest, making them unsuitable for the automatic synthesis of more complex, multi-dimensional geometries [17].

Local-search-based topology synthesis provides an efficient alternative to population-based methods. It involves adjusting one structure at a time, requiring only a limited number of EM simulations per iteration to identify appropriate descent direction for minimization of objective function. When embedded within a surrogate-assisted framework, the local-search process offers significantly reduced cost compared to more conventional algorithms. A common approach involves fine-tuning of topologies using a gradient-based routine embedded in a trust-region (TR) loop, which generates and optimizes a local approximation model around the best available solution [16], [19], [20]. Such a strategy has been successfully applied to the development of multi-parameter antennas, including cases involving variable-fidelity EM simulations and sequential dimensionality adjustments [16].

Despite proven usefulness of TR-based approaches, their effectiveness is limited to improving already promising initial designs rather than exploring the broader feasible space. These methods heavily rely on selection of an appropriate starting point, often derived from empirical equations, or known (hence, already existing) shapes, which constrains their topological flexibility. Additionally, due to complexity and multi-modal character of the problem, successful development of antennas is also a subject to availability of appropriate response-processing mechanisms. While local optimization—particularly when coupled with surrogate-assisted methods—can be cost-effective, it remains challenged by limitations of initial designs, topology selection, cost management, and convergence issues leaving the problem of unsupervised antenna design unresolved [11]-[21].

This work, which is an extension of a conference manuscript [1], explores two strategies for specification-oriented development of topologically agnostic planar antenna structures: (i) bandwidth-specific design and (ii) bandwidth-enhanced optimization. Both approaches classify quasi-randomly generated candidate designs, followed by their TR-based tuning within a variable-fidelity framework. The process involves optimization of topology represented in the form of the low-fidelity model followed by fine-tuning of the promising solutions using high-fidelity simulations. Frequency responses of the antenna structures are represented in the form of carefully selected feature points to ensure algorithm convergence. The considered design strategies have been demonstrated based on a total of three case studies concerning automatic development of planar antennas represented using 23, 33, and 53 independent parameters. The last one has been designed with respect to bandwidth-specific requirements, whereas generation of lower-dimensional structures is oriented towards maximization of bandwidth around a target frequency. Besides broadband operation, the optimized radiators exhibit unconventional radiation patterns. Owing to planar topology, the considered antennas might be utilized as components of in-door localization systems, or medical imaging devices (e.g., for breast cancer detection). A comparison of the considered design routines with other state-of-the-art approaches for developing patch-based radiators has also been conducted. Additionally, the feature-enhanced TR approach for bandwidth-specific design has been benchmarked against more conventional optimization methods.

## 2. Methodology

The proposed methodology for feature-assisted automatic development of antenna structures in a variable-fidelity setup involves the use of diverse tools that include electromagnetic (EM) simulation models (at two levels of accuracy), problem-specific representation of antenna responses, generation of suitable initial designs for optimization, as well as their surrogate-assisted optimization. The contents of this section include formulation of the design problem, explanation of the feature-based optimization concept, discussion of the optimization engine, as well as introduction of two strategies oriented towards development of antenna topologies for bandwidth-specific requirements and concerning enhancement of bandwidth. The section is concluded by a summary of the presented framework.

*2.1 Problem Formulation*

Let $R_f(x) = R_f(x, f)$ be the antenna response obtained for the high-fidelity EM simulation model over the specified frequency sweep $f$ for the given vector of input parameters $x$. The design problem can be formulated as a non-linear minimization task:

$$x^* = \arg\min_{x \in X} \left( U\left( \boldsymbol{R}_f(x) \right) \right) \quad (1)$$

where $x^*$ is the optimal design to be found, whereas $X$ denotes a feasible region of the search space; $U$ is a scalar objective function. Unfortunately, direct optimization of (1) is impractical due to high cost of $\boldsymbol{R}_f$ model simulations, as well as a large number of design parameters required to represent the topology-agnostic geometry. The cost can be reduced by embedding the process within a variable-fidelity setup where the high-fidelity model is substituted with its numerically cheap (yet less-accurate) low-fidelity counterpart $\boldsymbol{R}_c$. The latter one is characterized, e.g., by relaxed mesh density, as well as other simplifications [22]. The coarse EM model can be first used to narrow-down the search space to the region of interest, whereas fine-tuning of the geometry can be performed at the last stage of the design process using high-fidelity simulations. The cost of antenna design can be further reduced by replacing (1) with surrogate-assisted design procedure that gradually approximates the desired solution by solving:

$$x^{(j+1)} = \arg\min_{x \in X} \left( U\left( \boldsymbol{R}_s^{(j)}(x) \right) \right) \quad (2)$$

where $\boldsymbol{R}_s^{(j)}$, $j = 0,1, \ldots$, is the iteratively updated surrogate model constructed based on either low- or high-fidelity model simulations, whereas $x^{(j+1)}$ represents approximation of the $x^*$ design.

The problem considered in this work involves automatic generation of topology-agnostic planar antennas followed by their optimization in a variable-fidelity, surrogate-assisted setup. The latter is governed by a local trust-region-based algorithm that process antenna responses represented in the form of features [23], [24].

2.2  *Feature-Based Representation of Antenna Responses*

One of the main challenges pertinent to specification-oriented design of antenna structures represented using generic (shape-wise) simulation models involves intricate mutual relations between their input parameters (i.e., structure physical dimensions) on the performance characteristics. As a consequence, the radiator responses are highly non-linear (hence difficult to optimize) functions of frequency. The problem can be mitigated by shifting the design process into a feature-based domain where the key properties of the frequency response (from the optimization perspective) are represented using a set of carefully selected coordinate points [23]. In such a setup, the antenna properties are first extracted from the responses obtained at the given design (or set of designs) and then used for construction of an auxiliary model suitable for optimization [23], [25].

Let $F(x) = P(R(x))$ be the antenna performance characteristics expressed in terms of the feature points, where $P$ is the function used for their extraction. Note that, depending on the design stage, $R$ refers either to $R_c$, or $R_f$ model responses. The feature-based response $F = F(x)$ is defined as follows:

$$F = \begin{bmatrix} \boldsymbol{\omega} \\ \boldsymbol{S} \end{bmatrix} = \begin{bmatrix} \omega_1 & \cdots & \omega_m \\ S_1 & \cdots & S_m \end{bmatrix} \tag{3}$$

Each column in $F$ represents a pair of points that refer to frequency ($\boldsymbol{\omega} = \boldsymbol{\omega}(x)$) and level ($\boldsymbol{S} = \boldsymbol{S}(x)$) of the response. The pairs can be extracted according to the specific frequency point $\omega_m$, $m = 1, \ldots, M$, (e.g., pertinent to a local minimum or maximum of the reflection characteristic) and its corresponding level $S_m$. Alternatively, a pair of points can be obtained with respect to an appropriately defined level of response characteristic $S_n$ (e.g., at the edge of the bandwidth) and the frequency point associated with it $\omega_n$. (note that $n \neq m$, and $n \leq M$). Compared to the original response, its feature-based representation is a much less-nonlinear function of input parameters which simplifies the optimization problem [23], [25]. The meaning of frequency- and level-related coordinates, as well as the concept concerning representation of the response using features are illustrated in Figs. 1 and 2, respectively. The set of problem-specific points can be exploited to perform surrogate-assisted antenna tuning using the algorithm outlined below [23].

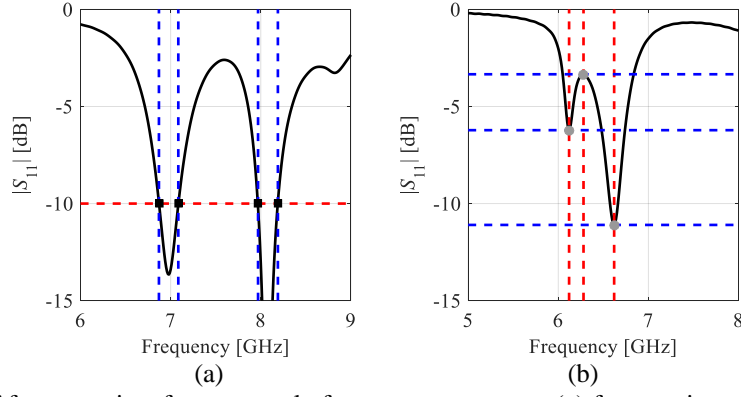

**Fig. 1.** Extraction of feature points from example frequency responses: (a) frequencies $\omega_{1-4}$ at pre-defined levels pertinent to reflection of $S_{1-4} = -10$ dB, as well as (b) response values $S_{1-2}$ and $S_3$ at frequencies $\omega_{1-2}$ and $\omega_3$ corresponding to minima and local maximum of the response. Blue lines denote the specific component of the response pertinent to the reference (red lines).

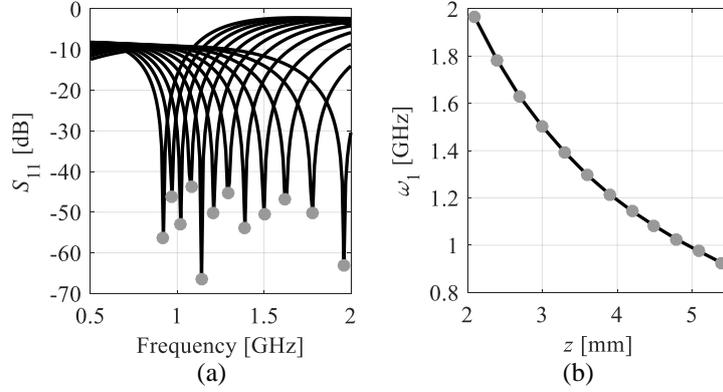

**Fig. 2.** Representation of structure response in the form of feature points: (a) family of frequency responses obtained as a function of some variable $x$, and (b) feature based response that corresponds to change of resonant frequencies with $x$. Note that variation of feature coordinates is much less non-linear compared to frequency characteristics.

### 2.3 Optimization Engine

The antenna is optimized using a gradient-based routine embedded in a trust-region (TR) framework. The algorithm implements a series of steps, $i = 0, 1, \ldots$, that involve: (i) preparation of the feature-based surrogate model, (ii) its optimization, (iii) verification of the approximated design solution, and (iv) re-set/refinement of the surrogate. Optimizations in the second step are performed by solving [20]:

$$\boldsymbol{x}^{(i+1)} = \arg \min_{\|\boldsymbol{x}-\boldsymbol{x}^{(i)}\| \leq \delta} \left( U\left(\boldsymbol{G}_\varepsilon^{(i)}(\boldsymbol{x})\right)\right) \tag{4}$$

The composite model $\boldsymbol{G}_\varepsilon^{(i)} = [\boldsymbol{G}_\omega^{(i)}\ \boldsymbol{G}_S^{(i)}]^T$, where $\boldsymbol{G}_\omega^{(i)} = \boldsymbol{\omega}(\boldsymbol{x}^{(i)}) + \boldsymbol{J}_\omega(\boldsymbol{x}^{(i)})(\boldsymbol{x} - \boldsymbol{x}^{(i)})$ and $\boldsymbol{G}_S^{(i)} = \boldsymbol{S}(\boldsymbol{x}^{(i)}) + \boldsymbol{J}_S(\boldsymbol{x}^{(i)})(\boldsymbol{x} - \boldsymbol{x}^{(i)})$ represent the first-order Taylor expansion surrogates generated around

$x^{(i)}$ w.r.t. the extracted frequency- and level-related features (3), respectively. The Jacobians are based on forward finite differences (FD) [23]:

$$J_\omega\left(x^{(i)}\right) = \begin{bmatrix} \cdots & \left(\omega\left(x^{(i)} + p_d^{(i)}\right) - \omega\left(x^{(i)}\right)\right)\frac{1}{p_d^{(i)}} & \cdots \end{bmatrix}^T$$
$$J_S\left(x^{(i)}\right) = \begin{bmatrix} \cdots & \left(S\left(x^{(i)} + p_d^{(i)}\right) - S\left(x^{(i)}\right)\right)\frac{1}{p_d^{(i)}} & \cdots \end{bmatrix}^T \quad (5)$$

Note that the parameter $p_d^{(i)}$ ($d = 1, \ldots, D$) denotes the FD perturbation w.r.t. $d$th dimension of the (currently best) design $x^{(i)}$, whereas for $p_d^{(i)}$ vector all components but the $d$th element—which is equal to $p_d^{(i)}$—are set to zero. The FD steps are relatively large so as to mitigate the effects of numerical noise on the quality of Jacobians while maintaining acceptable truncation errors [26]. This has been achieved by selecting the perturbations in proportion to the design $x^{(i)}$ and their update after each successful iteration [26]. The radius $\delta$ is controlled based on the ratio $\rho = [U(F(x^{(i+1)})) - U(F(x^{(i)}))]/[U(G_\varepsilon^{(i)}(x^{(i+1)})) - U(G_\varepsilon^{(i)}(x^{(i)}))]$ which is used to evaluate the obtained design solution and determine the procedure for surrogate update. The initial radius is set to $\delta = 1$. When $\rho < 0.25$ (poor prediction of $G_\varepsilon^{(i)}$) then $\delta = \delta/3$ which limits the expected range for which the surrogate model is considered appropriate, whereas for $\rho > 0.75$, $\delta = 2\delta$ (enhancement of range for which model is expected to be acceptably accurate) [23]. The gain coefficient $\rho$ is also used to accept ($\rho > 0$), or reject ($\rho \leq 0$) the candidate designs obtained from (4). The algorithm is terminated when $\delta^{(i+1)} < \varepsilon$, or $\|x^{(i+1)} - x^{(i)}\| < \varepsilon$ (here, $\varepsilon = 10^{-3}$).

It is worth emphasizing that the TR-based algorithm exploits surrogate models identified based on response features. In other words, the method works as follows. First, the frequency characteristics of the antenna are obtained around currently the best design $x^{(i)}$ and its FD-based perturbations. Then, the obtained data are used for extraction of feature coordinates $F(x^{(i)})$. The latter represents the input for construction of $G_\varepsilon^{(i)}$ model which undergoes optimization by solving (4). Next, the predicted design $x^{(i+1)}$ is evaluated based on EM simulation and the resulting frequency responses are processed to extract the new set of points $F(x^{(i+1)})$. Given that the antenna structures are topologically agnostic, the shapes at both designs and their corresponding frequency responses might be significantly different leading to difficulties in unequivocal identification of the number and location (frequency-, or level-wise) of relevant features (see Fig. 3). Here, the problem is mitigated by restricting extraction of $F(x)$ to a specific frequency range followed by identification of all pre-defined features confined within it (regardless of their number). At the same time, the objective functions used in the course of topology development are capable of handling inputs with variable length while identifying the performance-relevant coordinates. For more detailed discussion on TR-based optimization of problems represented in terms of response features, see [19], [20].

*2.4 Topology Generation and Optimization Strategies*

Determination of design optimization strategy is dictated by intended performance characteristics of the antenna at hand. Here, two strategies oriented towards generation of the topology-agnostic structures are considered. The first involves development of geometry dedicated to work within a specific, pre-determined bandwidth. The second is oriented towards maximization of the antenna bandwidth around the selected frequency of interest.

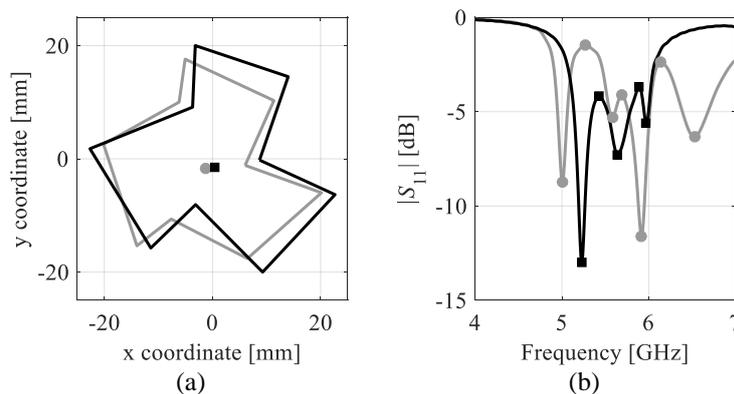

**Fig. 3.** Comparison of topology-agnostic antenna at two different iterations of TR-based optimization: (a) geometries of the radiator and (b) their corresponding frequency responses with highlights on extracted features (●) – here local minima and maxima over the frequency range of interest. Shape changes at the considered designs introduce substantial differences to the obtained performance characteristics. Consequently, the number of extracted features varies between the designs, i.e., five (black) vs. seven (gray). Such changes of frequency response properties between the iterations pose a practical challenge for unequivocal identification of feature points being of interest from the perspective of the design problem at hand.

### 2.4.1 Optimization for Bandwidth-Specific Applications

The first considered design strategy, originally published in [1], is a simplistic bi-stage procedure where the promising candidates are identified and then optimized using the algorithm of Section 2.3. The process is performed based on $R_c$ model simulations. The acceptable initial design is identified from a set of candidate solutions generated using a quasi-random method. The process involves determining a set of randomly distributed coordinates dedicated to represent topology-agnostic geometry (the generic EM model of the antenna used here re-creates the shape based on coordinate points), as well as their appropriate sorting in order to ensure shape convexity (and hence lack of its self-intersections). The obtained candidate solutions are then evaluated using a simple min-max classifier and sorted according to their performance.

The first step can be summarized as follows. Let $X_o = \{x_o\}$, $o = 1, \ldots, O$, $(X_o \subset X)$ be an $O$-element set of random designs that represent the geometrically-flexible antenna and $\{R_c(x_o)\}$ be the set of their corresponding low-fidelity EM-based responses. Each of the obtained designs can be evaluated as:

$$E_o = \max\left(R_c\left(x_o\right)\right)_{f_L \leq f \leq f_H} \tag{6}$$

where $f \in \boldsymbol{f}$ represent the frequency points within the bandwidth from $f_L$ to $f_H$ being of interest. Upon evaluation, the responses $E_1, \ldots, E_O$ are sorted according to their values and the design $x^{(0)} \in X_r$ that corresponds to the solution characterized by the lowest in-band reflection response is selected as a starting point for local TR-based optimization.

The bandwidth-specific optimization is handled using feature-specific objective function. Let $F(x) = [F_1(x)\ F_2(x)]$ where $F_1(x) = [\omega_1\ S_1]^T$ ($\omega_1 = [\omega_{1.1} \ldots \omega_{1.q}]$, $S_1 = [S_{1.1} \ldots S_{1.q}]$, $q = 1, \ldots, Q$) and $F_2(x) = [\omega_2\ S_2]^T$ ($\omega_2 = [\omega_{2.1} \ldots \omega_{2.b}]$, $S_2 = [S_{2.1} \ldots S_{2.b}]$, $b = 1, \ldots, B$; $Q < M$, $B < M$,) contain the feature points pertinent to local minima and maxima of the frequency response over the range of interest, respectively. The design objective is given as:

$$U(\boldsymbol{x}) = \sum \left( \max \left( S_2 - S_t, 0 \right)^2 \right) + \beta \| \boldsymbol{\omega}_1 - \boldsymbol{\omega}_t \| \tag{7}$$

where $S_t$ represents the threshold on reflection response and $\boldsymbol{\omega}_t = \omega_{t.1} + \delta \cdot [0 \ 1 \ \ldots \ Q{-}1]^T$ is a vector of target frequencies for the response minima. Given that the number of features might change in the course of TR-based optimization (cf. Section 2.3), the components of $\boldsymbol{\omega}_t$ represent evenly spaced target frequencies calculated based on the identified number of features $Q$ and the corner frequencies $\omega_{t.1} = f_L$ and $\omega_{t.Q} = f_H$ for the bandwidth of interest with $\delta = (\omega_{t.Q} - \omega_{t.1})/(Q - 1)$. The scaling coefficient $\beta = 100$ balances minimization of level- and frequency-related features for close-to-optimum solutions. Upon minimization of (7), the final obtained design $\boldsymbol{x}_c^*$ is used as a starting point for fine tuning. The process is also performed in a TR loop and governed using the function (7), yet based on the antenna responses evaluated using high-fidelity EM simulations.

### 2.4.2 *Optimization for Bandwidth-Enhancement*

The second design strategy introduces changes to identification of starting point for TR-based tuning, which are oriented towards increasing the amount of information on applicability of the design for optimization that can be extracted from the quasi-randomly generated solutions while maintaining acceptable cost of the process. The selected candidate designs undergo two step optimization which involves (i) shifting the antenna resonances towards the specific target frequencies and (ii) minimization of local maxima around the center frequency below the acceptable threshold while enhancing the bandwidth.

The modified topology generation procedure is based on the fact that antennas can be perceived as transformers between transmission lines and wireless propagation media [27]. Consequently, their electrical (and hence physical) dimensions are proportional to the wavelength propagated in the open space [27]. Furthermore, the probability that a randomly generated design will be characterized by suitable performance to be considered for optimization within a given frequency range of interest is relatively low (especially for multi-parameter geometries, as in this work). On the other hand, the geometry that (upon evaluation) does not produce resonances within the fraction of the frequency spectrum being of interest might feature promising responses (suggesting that the shape is useful for tuning) beyond it. The phenomenon can be leveraged through evaluation of the structure over a broad frequency range followed by uniform scaling of the radiator which would result in proportional shift of its response in frequency. Given floating-point nature of the approach, it is well supported by generic models represented in the form of coordinates (as in this work; cf. Section 3). Moreover, the scaling can be performed through interpolation of the antenna response according to an appropriately defined multiplicative factor. Regardless of inaccuracy compared to EM-based simulations of the scaled antenna model, the approach enables cost-efficient evaluation of the design candidate to the problem at hand.

Let $R_s(x_c, \alpha) = R_c(x_c, f_s)$ be an interpolated response of the topology-agnostic antenna obtained for the scaled frequency sweep $f_s = \alpha f$ at the quasi-randomly generated candidate design $x_c$. Then let $x^{(0)} = c \cdot x_c$ be the antenna geometry uniformly scaled using coefficient $c = 1/\alpha^*$ for which the response features the best performance w.r.t. the frequency range of interest. The optimum scale is obtained by solving:

$$\alpha^* = \arg\min_{\alpha} \left( U_c \left( R_s \left( x_c, \alpha \right) \right) \right) \tag{8}$$

The objective function $U_c = \max(R_s(x_c, \alpha)) - S_{max}$ is evaluated over a frequency range from $f_L$ to $f_H$ (cf. Section 2.4.1), whereas $S_{max}$ is a threshold that defines acceptable reflection of the scaled $x_c$ design within the band of interest. The candidate solution $x^{(0)}$ is accepted for optimization when the value of $U_c$ at $\alpha^*$ is at least equal to zero. Otherwise the design is discarded and a new candidate is generated for scaling. It is worth noting that evaluation of the candidate designs can be performed both sequentially, as well as in parallel (assuming availability of a suitable number of EM software licenses) so as to expedite the process of identifying a suitable starting point.

The feature-based optimization of the topology is performed as follows. Let $F(x) = [F_1(x)\ F_2(x)\ F_3(x)]$ represent the feature-based response of the antenna under design. The component $F_3(x) = [\omega_3\ S_3]^T$ ($\omega_3 = [\omega_{3.1} \ldots \omega_{3.g}]$, $S_3 = [S_{3.1} \ldots S_{3.b}]$, $g = 1, \ldots, G$) either denotes the frequency points at which the response passes through the $S_t$ threshold, or—if the $S_t$ is not crossed—is set to $F_3(x) = F_1(x)$. The meaning of $F_1(x)$ and $F_2(x)$ is the same as in Section 2.4.1. The design objective for the first step is $U_1 = \|\omega_1 - \omega_t\|$, where $\omega_t$ is obtained w.r.t. the specified corner frequencies and the number of identified features (cf. Section 2.4.1). The optimized design $x_{c.1}^*$ obtained as a result of $U_1$ minimization is then used as a starting point for the second stage of topology development. The objective function is given as:

$$U_2 = U_3 + \beta_1 \max\left(\max(S_2) - S_t, 0\right) - B \tag{9}$$

Here, $B = \min(B_1, B_2)$ with $B_1 = |f_0 - \omega_{3.1}|$, $B_2 = |\omega_{3.G} - f_0|$, whereas $f_0$ and $S_t$ denote the specified center frequency and the desired threshold on reflection response. The coefficient $\beta_1 = 10$ balances the level-related component of (9) against the frequency ones; $U_3$ is given as:

$$U_3 = \omega_{1.(Q-1)/2+1} \tag{10}$$

Similarly as for the first strategy, the design $x_c^* = x_{c.2}^*$ obtained as a result of feature-based optimization undergoes fine tuning based on high-fidelity EM simulations. The outline of variable-fidelity optimization using both considered strategies is given below.

*2.5 Variable-Fidelity Development of Antenna Designs*

The considered optimization framework combines quasi-random generation of topology-agnostic antenna designs and their optimization (both at the $R_c$ model level) using one of the selected strategies of Section 2.4. The final step involves re-optimization of the obtained low-fidelity designs in a TR setup based on high-fidelity EM simulations. As already mentioned, the goal of the considered variable-fidelity approach is to first explore the large search space pertinent to topology-agnostic design problems using computationally cheap simulations (yet characterized by certain inaccuracy) and then exploit the narrowed-down region so as to ensure that the design specifications are fulfilled. Given that the numerical cost of $R_f$ simulations is considerably higher compared to $R_c$, the approximation model $G_\varepsilon$ in (4) is constructed only once (at the beginning of optimization) and re-used upon termination of the optimization. Assuming that the feature extraction function $P$ is implemented (cf. Section 2.2), the algorithm can be summarized as follows (see Fig. 4 for conceptual illustration):

1. Set $D$, $f$, $f_L$, $f_H$, and $f_0$ (if needed);
2. Generate the candidate solution $x^{(0)}$ using the selected strategy (cf. Section 2.4);
3. Obtain $x_c^*$ through TR optimization of the $R_c$ model using the selected strategy;
4. Set $x^{(0)} = x_c^*$ and obtain $x_f^*$ using TR optimization of the $R_f$ model.

It should be noted that in Step 4 the design $x^{(0)}$ is already a good approximation of the final solution. Consequently, for the strategy of Section 2.4.2, the optimization involves only minimization of the objective function (9). For more comprehensive discussion on variable-fidelity optimization of antenna structures, see [28], [29], [30], [43].

## 3. Topology Agnostic Antenna

Automatic specification-driven antenna development is subject to availability of a generic model capable of supporting free-form geometries. Consider a topology-agnostic planar antenna shown in Fig. 5 [1]. The structure is implemented on a Rogers AD255C substrate with permittivity and thickness of $\varepsilon_r = 2.55$ and $h = 1.52$ mm, respectively. It comprises a patch fed through a concentric probe. The radiator is represented using the set of interconnected points defined in a cylindrical coordinate system. To ensure feasibility of geometry, each angular coordinate is set as relative to the previous one (the first is set to 0), whereas all angles are scaled w.r.t. their cumulative sum and normalized to the range between 0 and $2\pi$ (full circle). It should be noted that symmetry planes for the radiator are not defined so as to ensure its flexibility in terms of attainable performance specifications (even though this is achieved at the expense of increased computational cost associated with the size of EM model discretization domain) [22],

[30]. Note that such a representation supports quasi-random generation of topologies. For the probe, an additional routine is used to ensure that it is enclosed within the generated shape.

The antenna is implemented in a CST Microwave Studio and evaluated using its time-domain solver [31]. To ensure flexibility in terms of the number of design parameters, the structure EM model is generated dynamically and implements a safeguard mechanisms that verify consistency of the geometry and handle errors [31]. To support variable-fidelity design, two models of the structure have been implemented, i.e., $R_c$ (~100,000 tetrahedral mesh cells; simulation cost on an AMD EPYC 7282 system with 32 GB RAM: 60 s), and $R_f$ (~400,000 cells; simulation cost: 110 s). Contrary to $R_c$, the high-fidelity model implements lossy substrate as well as copper-based metallization with non-zero thickness.

The antenna geometry is represented using the following vector of design parameters $x = [C\ \rho_f\ \varphi_f\ \boldsymbol{\rho}\ \boldsymbol{\varphi}]^T$, where $C$ is a sizing factor for the patch/probe coordinates; $\rho_f$ and $\varphi_f$ represent the radial and angular position of the feed w.r.t. the origin of the cylindrical system, whereas the vectors $\boldsymbol{\rho} = [\rho_1\ ...\ \rho_l]^T$ and $\boldsymbol{\varphi} = [\varphi_1\ ...\ \varphi_l]^T$ ($l = 1, ..., L$) comprise the points that represent the shape of the patch (cf. Fig. 5). Note that the vector $x$ comprises $D = 2L + 3$ design parameters. The parameters $o = 5$, $r_1 = 1.27$, and $r_2 = 2.84$ remain fixed to ensure 50 Ohm input impedance, whereas $A = B = 2(C \cdot \max(\boldsymbol{\rho}) + o)$. All dimensions, except for the fixed ones and $C$ (in mm), are unit-less.

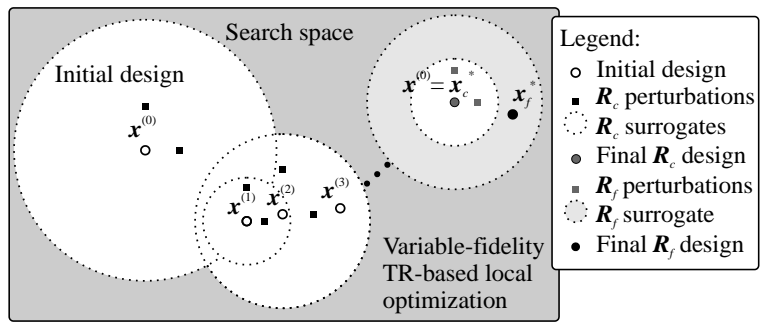

**Fig. 4.** Conceptual illustration of a variable-fidelity optimization using TR-based engine with linear model based on forward finite differences (cf. Section 2.3). The low-fidelity simulations are used in the course of search space exploration in order to narrow-down the region of interest (i.e., the one that contains optimum design). The identified fragment is then exploited using the surrogate constructed from $R_f$ model responses.

## 4. Results

In this section, the considered topology-agnostic antenna structure is optimized. A total of three case studies concerning generation of geometries for bandwidth-specific applications, as well as maximization of antenna bandwidth are discussed. The number of design parameters used in the course of structure development ranges from 23 to 53, respectively. The antenna topology is considered feasible within the appropriately specified lower and upper bounds: $l = [c_{0.l}\ 0\ \varphi_{0.l}\ 0.1\mathbf{1}\ 0.01\mathbf{1}]^T$ and $u = [c_{0.h}\ \rho_{0.f}\ \varphi_{0.h}\ 0.9\mathbf{1}\ 0.8\mathbf{1}]$T, where $\mathbf{1}$ is an $L$-dimensional vector of ones; $\varphi_{0.l} = \varphi^{(0)} - \pi$, $\varphi_{0.h} = \varphi^{(0)} + \pi$, and $\rho_{0.f} = \max(\boldsymbol{\rho}^{(0)})$ are the constraints determined based on the randomly generated design $x^{(0)}$ selected for numerical optimization based on the classification mechanisms discussed in Section 2.4. It should be noted that, when scaling of candidate designs is considered (cf. Section 2.4.2), the sizing factor of the antenna at $x^{(0)}$ is adjusted as $C' = c \cdot C$ to account for the frequency shift determined by minimization of (8). For the first case study, the bounds for scaling coefficient are set to $c_{0.l} = 25$ mm and $c_{0.h} = 35$ mm, respectively. For the bandwidth-enhancement design, the bounds are determined w.r.t. the scale of the quasi-randomly defined design as $c_{0.l} = C^{(0)} - 2$ and $c_{0.h} = C^{(0)} - 2$ (cf. Section 3). The discussion of the results, as well as comparison of structures against the antennas from the literature are given in Section 5.

### 4.1 Case 1 – Bandwidth-Specific Design

The first case study concerns development of the antenna topology using the strategy of Section 2.4.1. The corner frequencies for the optimization are set to $f_L = 6.2$ GHz and $f_H = 6.8$ GHz. The threshold on the reflection is set to $S_t = -10.2$ dB in order to provide a slight margin w.r.t. the intended in-band reflection level of $-10$ dB. The initial design $x^{(0)} = [30\ 0.13\ 2.62\ 0.22$ $0.36\ 0.48\ 0.57\ 0.59\ 0.56\ 0.53\ 0.44\ 0.34\ 0.36\ 0.35\ 0.43\ 0.52\ 0.54\ 0.38\ 0.29\ 0.43\ 0.42\ 0.42\ 0.47$ $0.5\ 0.57\ 0.41\ 0.29\ 0.22\ 0\ 0.17\ 0.18\ 0.26\ 0.26\ 0.26\ 0.29\ 0.25\ 0.3\ 0.29\ 0.42\ 0.3\ 0.31\ 0.05\ 0.05\ 0.31$ $0.13\ 0.43\ 0.36\ 0.31\ 0.29\ 0.04\ 0.14\ 0.3\ 0.6]^T$ is selected from a set of 200 randomly generated

candidate solutions evaluated using the min-max classifier (6). The optimized low-fidelity design $x_c^* = $ [30.04 0.13 2.47 0.21 0.36 0.47 0.56 0.61 0.58 0.52 0.46 0.34 0.36 0.31 0.44 0.53 0.55 0.38 0.33 0.44 0.41 0.42 0.49 0.52 0.57 0.41 0.28 0.23 0.01 0.16 0.17 0.25 0.26 0.27 0.29 0.24 0.3 0.29 0.42 0.29 0.29 0.08 0.07 0.32 0.13 0.43 0.36 0.31 0.29 0.04 0.15 0.32 0.61]$^T$ is found after 10 TR iterations. The final solution $x_f^* = $ [30.05 0.13 2.46 0.22 0.36 0.47 0.56 0.61 0.58 0.52 0.46 0.34 0.36 0.31 0.44 0.53 0.54 0.39 0.33 0.44 0.41 0.42 0.49 0.52 0.57 0.41 0.28 0.23 0.01 0.16 0.17 0.25 0.26 0.27 0.29 0.24 0.3 0.29 0.42 0.29 0.29 0.08 0.07 0.32 0.13 0.43 0.36 0.31 0.29 0.04 0.15 0.32 0.61]$^T$ has been obtained after 9 iterations of the algorithm. The computational cost of the design process corresponds to a total of 290.6 $R_f$ model simulations (~8.9 hours of CPU-time) including 200 $R_c$ and 221 $R_c$ evaluations for determination of initial design and its optimization, as well as 61 $R_f$ simulations for fine-tuning of the geometry.

Responses of the antenna at each stage of the design process, as well as visualization of the optimized design are shown in Fig. 6. The optimized geometry features –10 dB reflection within the operational range from 6.13 GHz to 6.83 GHz which provides a slight margin w.r.t. the intended bandwidth. Figure 7 demonstrates radiation patterns of the antenna obtained at $f_L$ and $f_H$. The structure features relatively high cross-polar gain at both considered frequencies which suggests its (local) close-to-circular polarization capabilities. Furthermore, the co-polar radiation around 6.8 GHz features a dual-lobe behavior with local minimum at 0° and maxima around ±45° angles, respectively. The maximum realized gain for the latter is around 5.11 dB at 40° and 2.26 at –50° for co- and cross-polarization. It is worth noting that the obtained unorthodox radiation pattern is only a by-product of the reflection-oriented optimization. Given that the structure has been developed to support 5$^{th}$ channel of ultra wideband frequency range, the radiation responses (upon further tuning) make it of potential use for radio-frequency-based in-door positioning. In particular, the radiator could be used for installation in convex corners within the facilities dedicated for in-door monitoring applications [9]. It is worth noting that, on a conceptual level, application of the structure in a heterogeneous localization architecture (with antennas tailored to conditions) could reduce the overall cost of UWB system installation [32].

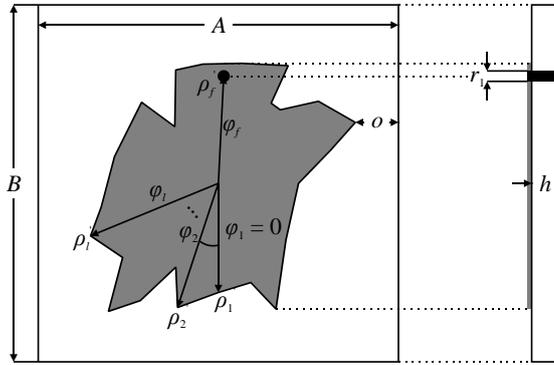

**Fig. 5.** Geometry of the proposed topologically agnostic patch antenna with highlight on the design parameters. Note that $\rho_f' = C\rho_f$ and $\rho_l' = C\rho_l$ ($l = 1, …, L$).

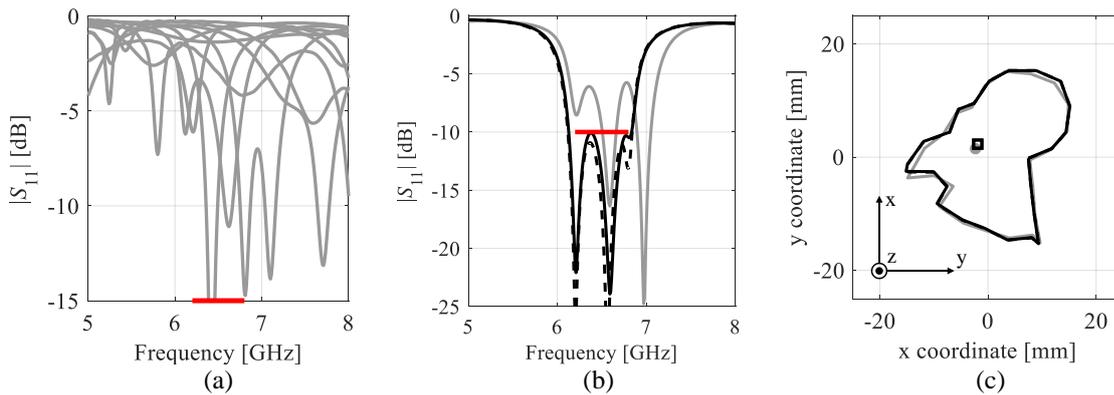

**Fig. 6.** Optimization of topology-agnostic antenna: (a) a few of randomly selected designs and bandwidth of interest for the min-max classifier (red line), (b) frequency responses at $x^{(0)}$ (gray), $x_c^*$ (– –) and $x_f^*$ (—), as well as (c) topology of the optimized structure.

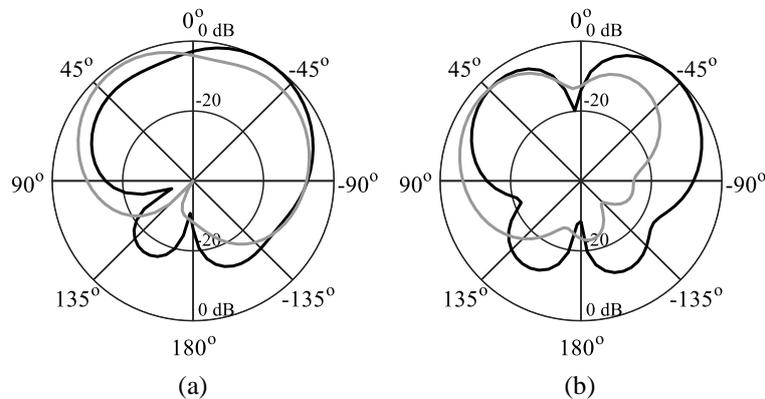

**Fig. 7.** Co- (black) and cross-polar (gray) radiation patterns of the topology-agnostic antenna obtained in the xz-plane (cf. Fig. 6) at: (a) 6.2 GHz and (b) 6.8 GHz frequencies.

*4.2 Case 2 – Bandwidth-Enhanced Optimization*

The second case study involves the development of an antenna topology with 23 design parameters using the strategy discussed in Section 2.4.2. The corner frequencies for the first stage of optimization process are set to $f_L$ = 5.35 GHz and $f_H$ = 5.65 GHz (with $f_0$ = 5.5 GHz). The threshold for the optimized reflection response is $S_t$ = –10.5 dB (again slightly lower compared to the intended in-band level of $S_t$ = –10 dB to account for modeling errors). The initial design $x^{(0)}$ = [44.90 0.03 3.66 0.14 0.36 0.43 0.28 0.45 0.44 0.27 0.40 0.47 0.14 0.00 0.78 1.04 0.24 0.90 0.99 0.11 1.05 0.95 0.23]$^T$ is obtained through interpolation of the design selected using the classifier with scaling. The optimized low-fidelity design $x_c^*$ = [44.76 0.04 3.89 0.13 0.43 0.47 0.23 0.54 0.46 0.18 0.45 0.51 0.15 0.04 0.75 0.80 0.26 0.80 0.77 0.11 0.79 0.79 0.25]$^T$ is found after 18 TR iterations using 320 $R_c$ evaluations. The final high-fidelity response $x_f^*$ = [44.76 0.04 3.89 0.13 0.43 0.47 0.23 0.54 0.46 0.18 0.45 0.51 0.13 0.00 0.88 0.94 0.31 0.93 0.91 0.13 0.93 0.92 0.29]$^T$ is obtained after just a few TR steps. The computational cost of the design process corresponds to a total of 284.3 $R_f$ model simulations (~8.7 hours of CPU-time) involving 93 $R_c$ and 320 $R_c$ evaluations for determination and optimization of the initial design, as well as a 58 $R_f$ simulations for fine-tuning the geometry.

The geometries and the reflection responses of the patch at the initial and optimized designs are compared in Fig. 8. The final geometry achieves a reflection of –10 dB within the operational range from 4.97 GHz to 5.72 GHz that represents a bandwidth of 0.75 GHz. It is worth noting that $x_c^*$ and $x_f^*$ are not centered w.r.t. $f_0$ due to lack of mechanism for shifting the local maximum in the case when even number of local minima is obtained through the optimization (as in this example). The far-field characteristics of the antenna that illustrate gain versus frequency at +35° and –32° angles (in the xz-plane) are shown in Fig. 9(a). The antenna consistently maintains a gain of at least 3 dB over the 5.1 GHz to 5.7 GHz frequency range in both directions. Additionally, the radiation patterns of the final optimized antenna in the co- and cross-polarization planes at 5.2 GHz and 5.6 GHz are shown in Figs. 9(b)-(c). The antenna exhibits a dual-lobe radiation pattern with a realized gain of 7.12 dB and 0.92 dB for the co- and cross-polarization patterns at 5.6 GHz, and corresponding values of 6.59 dB and –0.25 dB at 5.2 GHz, respectively. It is noteworthy that, similarly as in Section 4.1, the far-field parameters are merely by-products of the reflection-based optimization. The developed broadband antenna can be find applications for WiFi systems operating over the 5 GHz channels.

*4.3  Case 3 – Bandwidth-Enhanced Optimization*

The last example concerns application of the method outlined in Section 2.4.2 for development of the topology characterized by a total of 33 parameters. The first step of the optimization process targets corner frequencies of $f_L$ = 5.25 GHz and $f_H$ = 5.75 GHz (with reflection threshold set to $S_t$ = –10.5 dB – the same as for the antenna of Section 4.2). The initial design $x^{(0)}$ = [44.90 0.03 3.66 0.14 0.40 0.40 0.43 0.39 0.38 0.48 0.40 0.46 0.37 0.30 0.42 0.51 0.40 0.14 0.00 0.20 0.67 0.67 0.25 0.18 0.64 0.68 0.59 0.05 0.63 0.60 0.48 0.51 0.15]$^T$ is derived through interpolation applied to a geometry determined by the scaling-based classifier. The optimized low-fidelity design $x_c^*$ = [44.87 0.02 3.57 0.11 0.39 0.46 0.42 0.35 0.41 0.48 0.38 0.48 0.35 0.27 0.43 0.52 0.42 0.13 0.04 0.20 0.63 0.64 0.24 0.17 0.61 0.68 0.58 0.04 0.61 0.58 0.50 0.53 0.15]$^T$ is obtained after 18 TR iterations (a total of 449 $R_c$ evaluations). The final high-fidelity solution $x_f^*$ = [44.89 0.02 3.58 0.11 0.39 0.46 0.42 0.35 0.41 0.48 0.38 0.49 0.35 0.27 0.43 0.52 0.42 0.10 -0.00 0.19 0.64 0.65 0.24 0.17 0.62 0.68 0.59 0.04 0.62 0.58 0.51 0.54 0.16]$^T$ is found after 10 iterations of TR algorithm (110 $R_f$ evaluations). The computational effort for the design process amounts to 405.6 $R_f$ model simulations (~12.4 hours of CPU-time) involving 93 $R_c$ and 449 $R_c$ evaluations for the initial design and its optimization, along with 110 $R_f$ simulation for final geometry fine-tuning.

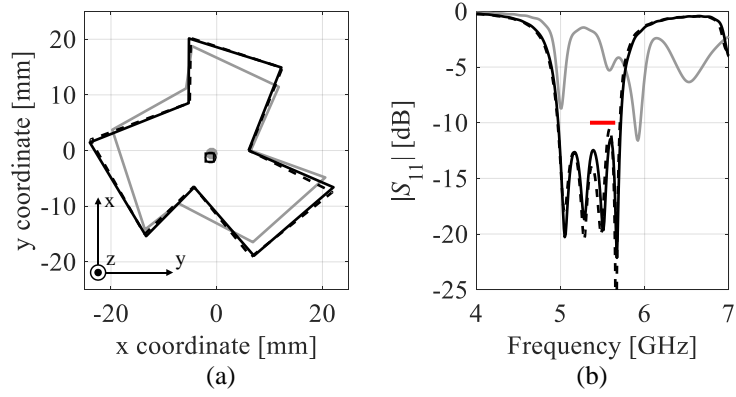

**Fig. 8.** Optimization of topology-agnostic antenna: (a) geometry, and (b) frequency responses at $x^{(0)}$ (gray), $x_c^*$ (– –), and $x_f^*$ (—), respectively.

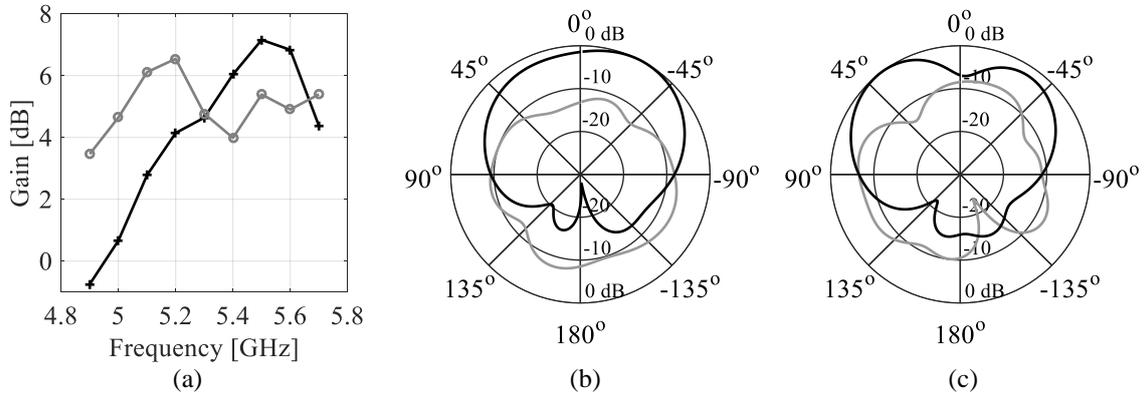

**Fig. 9.** The far-field responses for the antenna of Fig. 8 obtained in the xz-plane: (a) the gain versus frequency at directions of +35° (black) and –32° (gray), as well as radiation patterns in co- (black) and cross-polarization (gray) at: (b) 5.2 GHz and (c) 5.6 GHz, respectively.

Geometry and responses of the antenna structure at the initial and optimized designs are presented in Fig. 10. The final solution is characterized by a –10 dB reflection level over the frequency range from 5.24 GHz to 5.93 GHz, resulting in a bandwidth of 0.69 GHz. Figure 11 illustrates the far-field responses of the antenna over the band of interest, as well as at 5.5 GHz and 6 GHz frequencies, respectively. Similarly as for the remaining radiators, the structure exhibits a dual-lobe radiation pattern. The co- and cross-polarization gains are 7.57 dB and 0.52 dB at 5.5 GHz, as well as 6.25 dB and –0.41 dB at 6.0 GHz. Again, a broadband antenna with dual-lobe patterns could be advantageous in challenging environments where, e.g., omnidirectional antennas are insufficient to deliver the desired quality of service while maintaining relatively sparse (hence, cost-efficient) distribution.

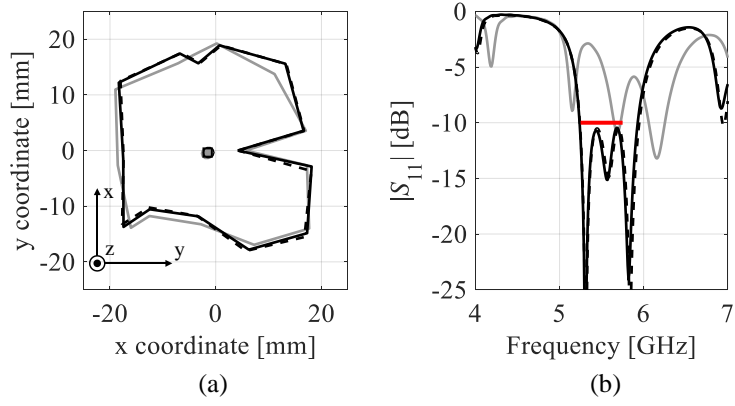

**Fig. 10.** Optimization of topology-agnostic antenna: (a) geometry, and (b) frequency responses at $x^{(0)}$ (gray), $x_c^*$ (– –) and $x_f^*$ (—).

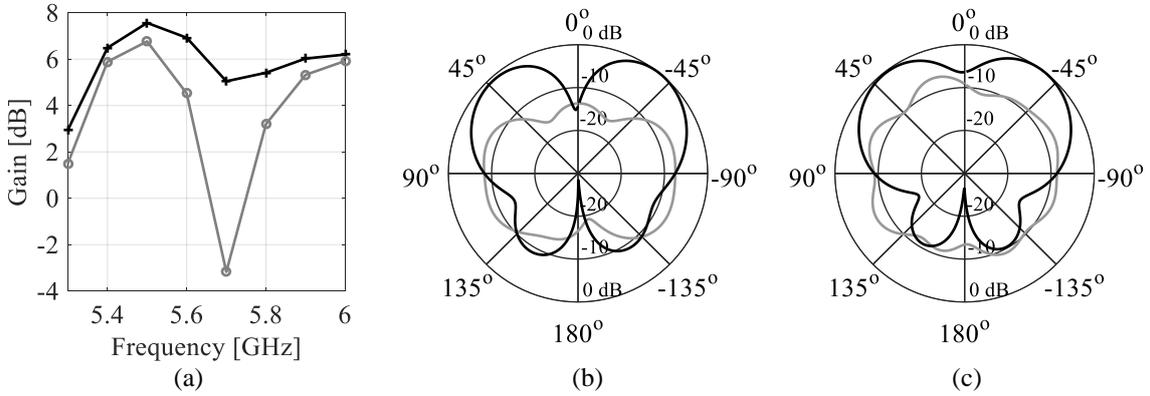

**Fig. 11.** Far-field responses for the antenna of Fig. 10 obtained in the xz-plane: (a) gain vs. frequency at +40° (black) and –40° (gray) angles, as well as co- (black) and cross-polarization (gray) radiation patterns at (b) 5.5 GHz and (c) 6.0 GHz, respectively.

## 5. Discussion and Comparisons

Generation of antenna topologies has been performed using two different strategies concerning optimization in feature-based setup for bandwidth-specific and bandwidth-enhanced operation. Regardless of different approaches and substantial differences between the resulting topologies, the generated antenna structures are characterized by similar electrical performance. The relative bandwidths for final designs from Sections 4.1 to 4.3 are 10.8% (0.7 GHz), 13.6% (0.75 GHz), and 12.5% (0.69 GHz), respectively. Furthermore, the results obtained in Section 4.2 indicate that the one of the driving factors behind enhanced bandwidth is increased number of local minima within the frequency range of interest rather than change of the number of design parameters. Such a conclusion can be drawn based on the observation that the structure with the lowest number of input variables features the broadest operational bandwidth. It should also be stressed out that the quasi-random generators utilized for both strategies do not guarantee that a suitable initial design for TR-based optimization will be found. On the other hand, it seems that the second classification strategy might outperform the one without uniform scaling in terms of identifying useful starting points for feature-assisted design. Clearly, the effect stems from shifting of local minima towards useful regions of the frequency spectrum.

Three developed antennas have been compared against the state-of-the-art radiators from literature in terms of size and performance [33]-[40]. The figures of interest include relative bandwidth, maximum gain, and electrical size. For fair comparison, the dimensions of each antenna have been expressed in terms of the guided wavelength $\lambda_g$, calculated with respect to the electrical parameters of the substrate underlying the radiator, as well as its center frequency $f_0$ [27]. The results gathered in Table 1 indicate that all but two of the considered antennas maintain narrower bandwidths than the generated topologies. Although the structure of [39] features low reflection for over 2 GHz (28%) along with relatively high gain, it has been implemented on a thick substrate and incorporates a multi-layer topology with several radiating components. Also, it is worth noting that the gain for antenna of [34] has not been reported. The radiators obtained in this work outperform most of the benchmark structures in terms of operational bandwidth and gain. However, the improved performance has been achieved at the expense of slightly larger electrical dimensions. Overall, the automatically developed radiators represent an acceptable compromise between physical dimensions and performance characteristics.

The discussed variable-fidelity optimization mechanism has been validated against competitive approaches from the literature. Specifically, the antenna of Section 4.1 has been optimized for minimization of the reflection within $f_L$ = 6.2 GHz to $f_H$ = 6.8 GHz bandwidth using the min-max function given as $U_0(x) = \max\{R(x)\}_{fL \leq f \leq fH}$, as well as the objective defined in a least-squares sense of the form:

$$U(x) = \frac{1}{N}\sum_{n=1}^{N} \max(R_n(x) - S_t, 0)^2 \qquad (11)$$

where $R_n \in R$ is the response of the antenna at the $n$th frequency point from $f_L$ to $f_H$ range. The considered benchmarks include direct optimization using $R_f$ model only in setups with: (i) min-max, (ii) least-squares, and (iii) feature-assisted objectives, as well as variable-fidelity approaches for (iv) min-max, and (v) function (11). For the sake of fair comparison, the starting point for each test is set to $x^{(0)}$ of Section 4.1. Furthermore, the cost concerning determination of the initial design is not accounted for in the comparison. The results gathered in Table 2 indicate that—due to complex topology of the antenna and intricate mutual relations between its design parameters—standard approaches that involve min-max or least-square optimization failed at obtaining acceptable design solution regardless of executing the process using high-, or low-fidelity EM-based simulations. At the same time, direct $R_f$-based tuning in a feature-based setup, while successful, is characterized by over 45% higher computational cost compared to optimization performed using the variable-fidelity setup. At the same time, both approaches produce radiators with similar performance in terms of in-band reflection, despite being characterized by slightly different topologies. A comparison of geometries and their responses obtained as a result of optimization using method (iii) and the one considered in Section 4.1 is shown in Fig. 12. It should be noted that both designs feature similar reflection performance despite being characterized by noticeably different topologies. The obtained results suggest that, for the considered design problem, feature-based optimization is a useful tool for determination of satisfactory antenna geometries.

Table 1: Comparison of the Optimized Designs

| Geometry | $f_0$ [GHz] | Bandwidth [GHz] | Gain [%] | [dB] | Footprint $[\lambda_g \times \lambda_g]$ | $[\lambda_g^2]$ |
|---|---|---|---|---|---|---|
| [33] | 5.90 | 0.27 | 4.58 | 6.02 | 0.40 × 0.60 | 0.24 |
| [34] | 5.50 | 0.87 | 15.8 | – | 0.46 × 0.65 | 0.30 |
| [35] | 2.46 | 0.18 | 7.30 | 2.90 | 0.55 × 0.57 | 0.31 |
| [36] | 6.70 | 0.35 | 5.23 | 1.54 | 0.49 × 0.65 | 0.32 |
| [37] | 3.86 | 0.08 | 2.17 | 5.34 | 0.59 × 0.85 | 0.50 |
| [38] | 5.30 | 0.31 | 5.83 | 1.48 | 0.93 × 0.93 | 0.87 |
| [39] | 7.52 | 2.11 | 28.0 | 8.00 | 0.79 × 2.08 | 1.64 |
| [40] | 5.40 | 0.17 | 3.05 | 7.90 | 1.36 × 1.36 | 1.85 |
| Section 4.1 | 5.50 | 0.70 | 10.80 | 5.11 | 1.10 × 1.02 | 1.12 |
| Section 4.2 | 5.20 | 0.75 | 14.42 | 6.59 | 1.21 × 1.21 | 1.46 |
| Section 4.3 | 5.60 | 0.69 | 12.32 | 7.02 | 1.03 × 1.03 | 1.06 |

Table 2: Benchmark of the Algorithm

| Method | $R_c$ | $R_f$ | Total cost $R_f$ | $R_f$ [h] | $U_0(x)$ [dB] |
|---|---|---|---|---|---|
| (i) | – | 61 | 61 | 1.86 | –5.58 |
| (ii) | – | 222 | 222 | 6.78 | –8.24 |
| (iii) | – | 331 | 331 | 10.1 | –10.3 |
| (iv)* | 275 | – | 150 | 4.58 | –6.93 |
| (v)* | 385 | – | 210 | 6.42 | –8.64 |
| Section 4.1 | 221 | 61 | 181.5 | 5.55 | –10.1 |

* Optimized design violates the specification – $R_f$ tuning not performed

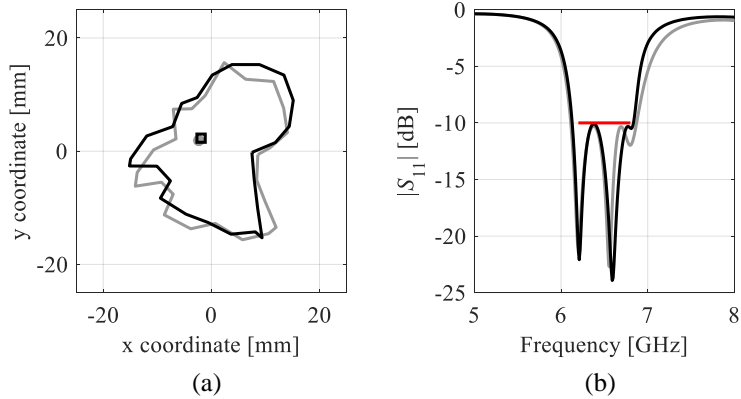

(a)      (b)

**Fig. 12.** Benchmark of the design algorithms – comparison of the optimized antenna designs obtained using the $R_f$-based feature-assisted method (gray) and the variable-fidelity algorithm of Section 2 (black) in terms of: (i) topology shape and (b) reflection responses at the final designs.

## 6. Conclusions

In this work, two strategies for feature-assisted development of topology-agnostic planar antennas have been considered. Both methods involve quasi-random generation of antenna topologies followed by their evaluation using appropriate classification functions and feature-based optimization using a trust-region method. The design method is embedded in a variable-fidelity framework where the topology is first optimized using low-fidelity model with relaxed mesh density and then fine-tuned based on high-fidelity model simulations. The considered design strategies have been validated based on a total of three case studies concerning development of a point-based antenna w.r.t. bandwidth-specific requirements and two generic geometries for maximization of bandwidth around the center frequency of interest. The average computational cost of the design using considered methods amounts to 327 $R_f$ simulations (~10 hours of CPU-time). It is worth noting that the optimized antenna designs feature diversified radiation patterns including the ones with a dual-lobe behavior at the selected frequencies. Such properties (upon suitable tuning) make them of potential use for in-door positioning in challenging propagation conditions. The results of comparisons indicate that the optimized structures often outperform state-of-the-art patch-based radiators from the literature in terms of bandwidth while maintaining acceptable dimensions. Furthermore, the strategy concerning feature-enhanced TR design of the antenna for bandwidth-specific applications has been favorably benchmarked against more conventional optimization methods. Future work will focus on extension of bandwidth-enhanced-oriented method so as to account for centering the response w.r.t. the center frequency regardless of the number of local minima.

## Acknowledgement

This work was supported in part by the National Science Center of Poland Grant 2021/43/B/ST7/01856.